\documentclass[a4paper,8pt]{article}

\usepackage[english]{babel}

\usepackage{natbib}

\usepackage{amsfonts,amsbsy,amssymb,amsmath,graphicx,float} 

\usepackage[paperwidth=18cm,paperheight=29.7cm,hmargin=2cm,vmargin=3cm]{geometry}

\usepackage[linktoc=all,colorlinks=true,citecolor=blue]{hyperref}

\begin{document}

\title{Unveiling the Fractal Structure \\ of Julia Sets with Lagrangian Descriptors}

\author{V\'{i}ctor J. Garc\'{i}a-Garrido$^{1}$ \\[.3cm]
	$^1$ Departamento de F\'{i}sica y Matem\'{a}ticas, \\ Universidad de Alcal\'{a}, 28871, Alcal\'{a} de Henares, Spain. \\[.2cm] vjose.garcia@uah.es}

\maketitle

\begin{abstract}
	
	In this paper we explore by means of the method of Lagrangian descriptors the Julia sets  arising from complex maps, and we analyze their underlying dynamics. In particular, we take a look at two classical examples: the quadratic mapping $z_{n+1} = z_n^2 + c$, and the maps generated by applying Newton's method to find the roots of complex polynomials. To achieve this goal, we provide an extension of this scalar diagnostic that is capable of revealing the phase space of open maps in the complex plane, allowing us to avoid potential issues of orbits escaping to infinity at an increasing rate. The simple idea is to compute the $p$-norm version of Lagrangian descriptors, not for the points on the complex plane, but for their projections on the Riemann sphere in the extended complex plane. We demonstrate with several examples that this technique successfully reveals the rich and intricate dynamical features of Julia sets and their fractal structure.  
	
\end{abstract}


\section{Introduction}
\label{sec:intro}

Complex dynamics \cite{milnor2011,carleson2013,beardon2000} is a branch of Dynamical Systems Theory concerned with the study of discrete-time maps generated by the iterative application of a complex function in one or several variables. This subject, which began in the first quarter of the twentieth century with the pioneering works by Pierre Fatou and Gaston Julia about the iteration of rational functions, gave birth to fractal geometry and led to the magnificent discovery of the most representative icon of infinite complexity, the Mandelbrot set, a true mathematical gem of incalculable beauty. For a fascinating account on the history and developments of this discipline, see \cite{alexander2013}. 

In this context, one of the biggest challenges of Nonlinear Dynamics is the development of mathematical techniques that provide us with the capability of exploring the geometrical template of structures that distinguishes regions with qualitatively distinct dynamical behavior in phase space. In the particular case of complex dynamics, the phase space is identified with the comple plane itself. Since the early 1900, the idea of pursuing a qualitative description of the solutions of differential equations that emerged from the work carried out by Henri Poincar\'e on the three body problem of Celestial Mechanics \cite{hp1890}, has had a profound impact on our understanding of the nonlinear character of natural phenomena.

The goal of this paper is to present the details of a mathematical scalar diagnostic for the analysis of dynamics in the complex plane. This method is known in the literature as Lagrangian descriptors (LDs) and it has been widely applied to problems that arise in many scientific disciplines. This tool was originally developed for the study continuous-time dynamical systems with general time-dependece \cite{madrid2009,mancho2013lagrangian,lopesino2017}. In this setup, it has been used for instance in Oceanography for the assessment of marine oil spills \citep{gg2016}, and also to plan autonomous underwater vehicle transoceanic missions \cite{ramos2018}. The versatility provided by this approach has also been evidenced in the adaption of this technique to explore the phase space of stochastic dynamical systems \cite{balibrea2016lagrangian}.

Lagrangian descriptors have also been defined in the setting where the dynamics of a system takes place in discrete time steps (iterations of a map). This work was first carried out in \citep{carlos}, and Discrete Lagrangian Descriptors (DLDs) have provided since then insightful results when used to describe 2D chaotic maps such as the H\'{e}non, Lozi, and Arnold's cat maps \citep{carlos,carlos2,GG2018RCD2}. In all these exaamples, that display strong mixing properties, DLDs have successfully recovered the chaotic saddles and strange attractors responsible for the chaotic dynamics. 

However, when Lagrangian Descriptors (both in their continous- or discrete-time versions) are used to explore the phase space structures of unbounded dynamical systems, issues might arise when trajectories of initial conditions escape to infinity in finite time or at an increasing rate. This behavior will result in NaN values in the LD scalar field, obscuring the detection of invariant manifolds. Recent studies \cite{junginger2017chemical,naik2019b,GG2019} have revealed this concern for some examples of unbounded continuous-time dynamical systems, where the approach of calculating LDs by integrating initial conditions on a phase space grid for the same time, known as fixed-time LDs, does not work well because of escaping trajectories. In order to circumvent this issue, an alternative definition of LDs was introduced for continuous-time systems, where the value of LDs for a given initial condition is calculated until the trajectory leaves a sufficiently large domain in phase space. This idea has also been implemented for the analysis of discrete-time maps in \cite{GG2019b} and it has shown to produce excellent results for unbounded maps. It is interesting to point out that an alternative extension of DLDs to study semiclassical mechanics and quantum chaos has also been given in \cite{carlo2019}. 

As we have mentioned, the objective we pursue in this work is to give an alternative definition of DLDs with the capability of revealing the internal dynamics of complex maps. In this task, we could have decided to directly apply the idea introduced in \cite{GG2019b} to stop the iterations of orbits when they leave a fixed phase space region, so that those that escape to infinity very fast do not ruin the interpretation of the output obtained from the scalar diagnostic. On the contrary, we have opted for a different strategy that relies on the construction of the extended complex plane and the computation of $p$-norm distances between points on the Riemann sphere. 

This paper is outlined as follows: Section \ref{sec:sec1} starts by briefly explaining the construction of the extended complex plane by means of applying the stereographic projection to the Riemann sphere. This is a standard procedure followed in complex dynamics that is suitable for the rigorous mathematical analysis of the convergence of sequences obtained from the iteration of maps in the complex plane. Once this framework has been established, we introduce an alternative definition of DLDs
that accumulates the $p$-norm along the orbits of initial conditions, and it does so not for the points on the complex plane, but for their projections on the Riemann sphere. Next, in Section \ref{sec:res} we discuss the results of this work. First, we numerically show how the adapted DLD is capable of revealing the geometrical template of phase space structures that determine the dynamics of the classical quadratic map generated by the complex function $f(z) = z^2 + c$. In order to provide further validation, and as a demonstration of the potential that this technique could bring, we apply it to unveil the fractal dynanmics induced by the use of Newton's method for the computation of the roots of complex polynomials. Finally, in Section \ref{sec:conc} we present the conclusions of this work.

\section{Discrete Lagrangian Descriptors for Complex Maps}
\label{sec:sec1}

In this section we construct a Discrete Lagrangian Descriptor (DLD) to explore the dynamics of complex maps in the form:
\begin{equation}
z_{n+1} = f(z_{n}) \quad,\;\; \forall \; n \in \mathbb{N} \cup \lbrace0\rbrace \quad , \;\; \text{where} \;\; z_n = x_n + y_n \, i \in \mathbb{C} \;.
\label{discrete_DS}
\end{equation}
One of the most important issues that must be addressed when analyzing the dynamics of such discrete time systems is that orbits of initial conditions can easily escape to infinity very quickly, making the study of the underlying dynamical structures a delicate task. Our goal is to introduce a scalar diagnostic that sums the values taken by a positive function along trajectories of initial conditions, and whose output reveals the geometrical structures in the phase space that determine the dynamics of the system.

The first step that we take in order to avoid issues of trajectories escaping to infinity is to compactify the complex plane by wrapping it around a sphere, and identifying all the points at infinity with the north pole of the sphere. This is done by using a stereographic projection that maps the points on the surface of a unit sphere, known as the Riemann sphere, onto its equatorial plane that is identified with the complex plane. The result of this construction is the extended complex plane, denoted by $\hat{\mathbb{C}} = \mathbb{C} \cup \left\{\infty\right\}$, and the Riemann sphere is its geometrical representation.

Consider a sphere with unit radius centered at the origin, and take a point $P$ on its surface with coordinates $\mathbf{\xi} = (\xi_1,\xi_2,\xi_3)$. The Riemann sphere is defined by:
\begin{equation}
\mathcal{R} = \left\{\mathbf{\xi} = (\xi_1,\xi_2,\xi_3) \in \mathbb{R}^3 \; \big| \; \xi_1^2 + \xi_2^2 + \xi_3^2 = 1 \right\}
\end{equation}
The equatorial plane of this sphere, i.e. $\xi_3 = 0$, represents the complex plane and any complex number $z = x + y \,  i \in \mathbb{C}$ can be seen as a point $P^{\prime}$ with cartesian coordinates $(x,y)$, where the $x$ and $y$ axes are taken along the $\xi_1$ and $\xi_2$ axes respectively. The stereographic projection of any point $P$ on the Riemann sphere onto the complex plane is easily obtained by constructing the straight line that joins the north pole of the sphere with $P$, and determining its intersection with the equatorial plane. This operation yields the following mapping:
\begin{equation}
\begin{array}{ccccc}
\mathcal{S} & : & \mathcal{R} & \quad \longrightarrow & \mathbb{C} \\[.2cm]
& & \xi = (\xi_1,\xi_2,\xi_3) & \quad \longmapsto & \quad z = \mathcal{S}(\xi) =  \dfrac{\xi_1 + \xi_2 \, i}{1-\xi_3}
\end{array}
\label{stp_map}
\end{equation}
with an inverse function that takes points from the complex plane to the Riemann sphere and is given by:
\begin{equation}
\begin{array}{ccccc}
\mathcal{S}^{-1} & : & \mathbb{C} & \quad \longrightarrow & \mathcal{R} \\[.2cm]
& & z = x + y \, i & \quad \longmapsto & \quad \xi = \mathcal{S}^{-1}(z) = \dfrac{1}{|z|^2 + 1}\left(2x,2y,|z|^2 - 1\right)
\end{array}
\label{stpinv_map}
\end{equation}

The method of Discrete Lagrangian Descriptors was originally introduced in \cite{carlos} in order to analyze the invariant manifolds that characterize the dynamics in the phase space of two-dimensional invertible maps. It was defined as follows: given a fixed number of iterations $N > 0$ of a map $f$ (forward and backward), and let $p \in (0,1]$ be the order that specifies the $l^p$-norm used to measure the distance between sequential iterations of the mapping. Then, the DLD is given by the function:
\begin{equation}
MD_{p}\left(\mathbf{x}_0,N\right) = \sum_{k=-N}^{N-1} ||\textbf{x}_{k+1}-\textbf{x}_{k}||_{p} = \sum_{k=-N}^{N-1} |x_{k+1}-x_{k}|^{p} + |y_{k+1}-y_{k}|^{p}
\label{DLD}
\end{equation}
where $\mathbf{x}_0 \in D$ is any initial condition chosen on a bounded subset $D$ of the plane, and $\mathbf{x}_k = f^{(k)}(\mathbf{x}_0)$ denotes the $k$-th iteration of the map. Observe that Eq. (\ref{DLD}) can be split into two quantities,
\begin{equation}
MD_{p}\left(\mathbf{x}_0,N\right) = MD_{p}^{+}\left(\mathbf{x}_0,N\right) + MD_{p}^{-}\left(\mathbf{x}_0,N\right) \;,
\end{equation}
where,
\begin{equation}
MD_{p}^{+} = \sum_{k=0}^{N-1} ||\textbf{x}_{k+1}-\textbf{x}_{k}||_{p} = \sum_{k=0}^{N-1} |x_{k+1}-x_{k}|^{p} + |y_{k+1}-y_{k}|^{p},
\label{DLD+}
\end{equation}
measures the forward evolution of the orbit of $\mathbf{x}_0$ and,
\begin{equation}
MD_{p}^{-} = \sum_{k=-N}^{-1} ||\textbf{x}_{k+1}-\textbf{x}_{k}||_{p} = \sum_{k=-N}^{-1} |x_{k+1}-x_{k}|^{p} + |y_{k+1}-y_{k}|^{p},
\label{DLD-}
\end{equation}
accounts for the backward evolution of the orbit of $\mathbf{x}_0$.
Notice that in this definition, the method takes into account the fact that the function defining the discrete dynamical system is invertible. However, in the complex mappings that we consider in this work, the functions that define the dynamics are not invertible, and therefore we will only keep the forward part of the definition. It is important to remark here that in \cite{carlos} it was mathematically proved that DLDs highlight stable and unstable manifolds at points where the $MD_p$ field becomes non-differentiable (its gradient becomes discontinuous or unbounded). These non-differentiable points are displayed in a simple way when depicting the $MD_p$ scalar field, and are known in the literature as ''singular features''. Therefore, DLDs are capable of recovering the geometrical template of invariant stable and unstable manifolds present in the phase space of any map. Moreover, \cite{carlos} shows for several examples that $MD_{p}^{+}$ detects the stable manifolds of the fixed points, whereas $MD_{p}^{-}$ does the same for the unstable manifolds of the fixed points. Therefore, the fixed points with respect to the map can be located at the intersections of these structures.

Another relevant aspect to remark from the definition of DLDs is that the number of iterations $N$ used to compute the orbit of any initial condition plays a crucial role for revealing the intricate structure of stable and unstable 
manifolds in the phase space. We can claim that the parameter $N$ holds the key to unlock the general applicability of DLDs to unbounded maps. In this respect, it has been also shown in \cite{carlos} that, as the number of iterations $N$ is increased, 
a richer template of geometrical phase space structures is recovered by the method, since we are including into the analysis more information about the past and future history of the orbits of the map. An issue could arise when one deals with an unbounded map for which orbits can escape to infinity at an increasing rate. With the purpose of circumventing this difficulty, DLDs were recently extended in \cite{GG2019b} in order to improve its capability for analyzing unbounded maps. This goal was achieved by accumulating the value of DLDs along the orbit of each initial condition up to a fixed number of iterations, or until the trajectory leaves a certain fixed region in the plane. For further details on how this modification was implemented for the method, see \cite{GG2019b}.

At this point we are ready to introduce a Discrete Lagrangian Descriptor specifically adapted for the analysis of the dynamics of complex map with the general form described in Eq. \eqref{discrete_DS}. Given a fixed number of iterations $N$ and a value $p \in (0,1]$ that determines the $l^p$-norm used to compute distances, we define the DLD as:
\begin{equation}
\mathcal{D}_{p}^{+}\left(z_0,N\right)
= \sum_{k=0}^{N-1} \big|\big|\mathcal{S}^{-1}(z_{k+1}) - \mathcal{S}^{-1}(z_{k})\big|\big|_p
= \sum_{j=1}^{3} \sum_{k=0}^{N-1} \big|\xi_j^{(k+1)}-\xi_j^{(k)}\big|^{p},
\label{DLD+complexmap}
\end{equation}
where $z_0 = x_0 + y_0 \, i \in D$ is any initial condition selected on a bounded subset $D$ of the complex plane and $\xi^{(k)} = (\xi_1^{(k)},\xi_2^{(k)},\xi_3^{(k)}) = \left(S^{-1} \circ f^{(k)}\right)(z_0)$ represents the projection of the $k$-th iterate of the map onto the Riemann sphere. It is important to observe that, since the inverse of the stereographic projection defined in Eq. \eqref{stpinv_map} is a smooth function, the singular features highlighted by the method at the points were the scalar field $MD_{p}$ is non-differentiable wil be preserved. This ensures that the alternative definition that we are giving still recovers the invariant manifolds in the phase space of the complex map.

\section{Results}
\label{sec:res}

We begin our discussion of the results by applying LDs to analyze the complex dynamics of the classical quadratic map $z_{n+1} = z_n^2 + c$, where $c = \alpha + \beta \, i \in \mathbb{C}$ is a complex parameter. This yields the corresponding two-dimensional discrete dynamical system in the plane:
\begin{equation}
\begin{cases}
x_{n+1} = x_n^2 - y_n^2 + \alpha \\[.2cm]
y_{n+1} = 2 x_n \, y_n + \beta
\end{cases}
\, , \quad z_n = x_n + y_n \, i \in \mathbb{C} \;.
\label{ds_qwad}
\end{equation}
Before analyzing the complex dynamics of this map, we list all the parameter values that we look at:
\begin{itemize}
	\item $c = 0$ (Unit Circle Julia Set).
	\item $c = i$ (Dendrite Fractal).
	\item $c = 0.25$ (Cauliflower Parabolic Set).
	\item $c = -1$ and $c = -0.75$ (San Marco Basilica-type Julia Sets).
	\item $c = -0.123 + 0.745 \, i$ (Douady's Rabbit).
	\item $c = 0.285 + 0.01 \, i$ (Open Culiflower Set).
	\item $c = -0.1 - 0.651 \, i$ (Flower-type Fractal).
	\item $c = -0.391 - 0.587 \, i$ (Siegel Disk Fractal).
\end{itemize}

\begin{figure}[!h]
	\begin{center}
		A)\includegraphics[scale=0.31]{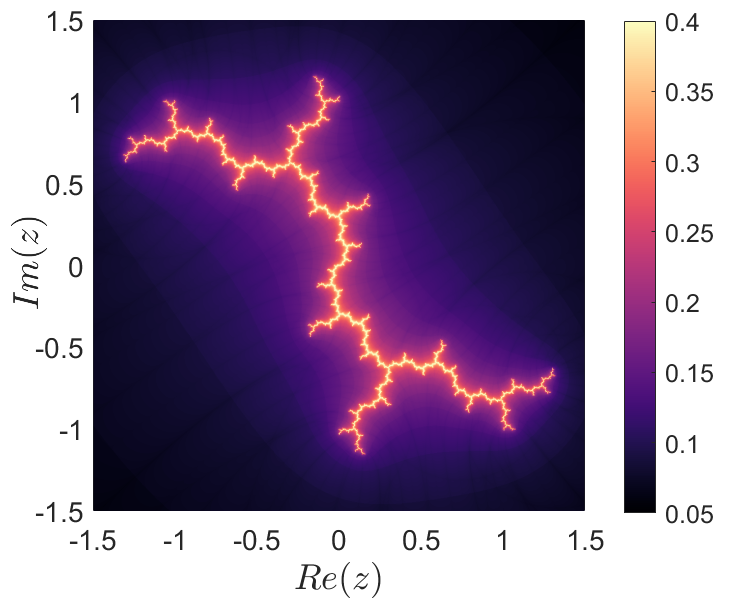}
		B)\includegraphics[scale=0.46]{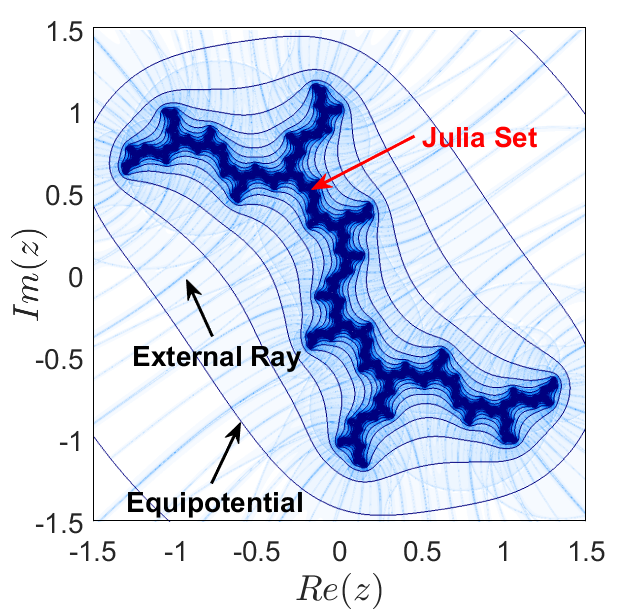}
	\end{center}
	\caption{A) Dendrite fractal obtained for the quadratic mapping with $c = i$, as revealed by DLDs using $p = 1/64$ and $N = 200$ iterations. B) Julia fractal set, equipotential curves and external rays extracted from the gradient of Lagrangian descriptors.}
	\label{LDs_dendrite}
\end{figure}

Consider first the case where $c = i$. For this value, the Julia set obtained is known as a dendrite, which are Julia sets with no interior that display a tree-like branching structure. In order to reveal this structure we compute DLDs using $p = 1/64$ for $N = 200$ iterations. The output of the method is shown in Fig. \ref{LDs_dendrite}, where the scalar field depicted in A) nicely captures the fractal set. One of the most important features of DLDs is that one can extract the key dynamical features present in the phase space of the map from the ridges of the gradient $||\nabla \mathcal{D}_{p}^{+}||$. This is so, because in \cite{carlos} it was mathematically proven that when DLDs are applied in forward time, the stable manifolds are revealed at points where the function is non-differentiable, and this implies that the gradient becomes discontinuous or unbounded. Following this approach, we calculate in Fig. \ref{LDs_dendrite} B) the points where the gradient attains very large values, and this procedure gives us, not only the Julia set corresponding to the dendrite, but also the equipotential lines and external rays of the complex mapping.  

\begin{figure}[!h]
	\begin{center}
		A)\includegraphics[scale=0.46]{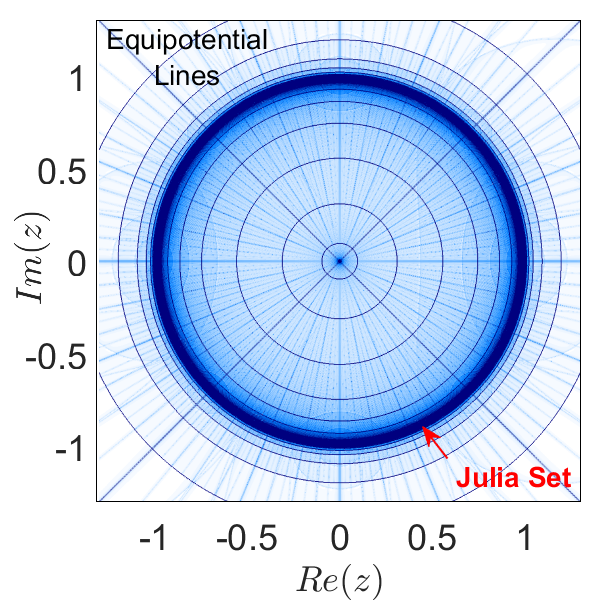}
		B)\includegraphics[scale=0.54]{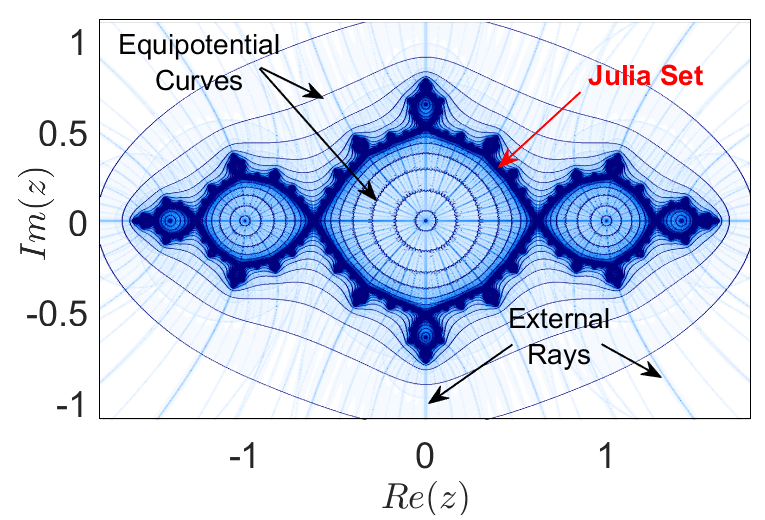}
	\end{center}
	\caption{Phase space of the quadratic map revealed by DLDs with $p = 1/64$ and $N = 1000$ iterations for: A) $c = 0$; B) $c = -1$.}
	\label{LDs_equip_extrays}
\end{figure}

\noindent With the purpose of a further validation that the method is faithfully recovering relevant structures for the analysis of complex dynamics, we apply DLDs to the quadratic mapping with $c = 0$ and $c = -1$. The case $c = 0$, despite being very simple from the point of view of the phase portrait, it displays chaotic dynamics. The Julia set associated to this parameter value is the unit circle. All initial condtions selected on it can be written as $z_0 = e^{i \theta}$, and those with an irrational argument behave chaotically and their corresonding orbit fills up the entire circle uniformly. On the other hand, points inside the unit circle give rise to sequences that converge towards the origin, which is a stable fixed point of the map. Notice that the fixed points are revealed at points for which the DLD is zero, since they are stationary. The points located outside the unit circle tend to infinity, that is also a fixed point. Observe that the method gives us the equipotential curves as well when we depict the ridges of large gradient values. All this analysis is displayed in Fig. \ref{LDs_equip_extrays} A). The computaion of DLDs using $p = 1/64$ and $N = 1000$ iteartions for the case $c = -1$ reveals in Fig. \ref{LDs_equip_extrays} B) the San Marco fractal, its equipotentials and also the external rays \cite{belk2015}. This Julia set takes its name from San Marco Basilica in Venice.

\begin{figure}[!h]
	\begin{center}
		A)\includegraphics[scale=0.28]{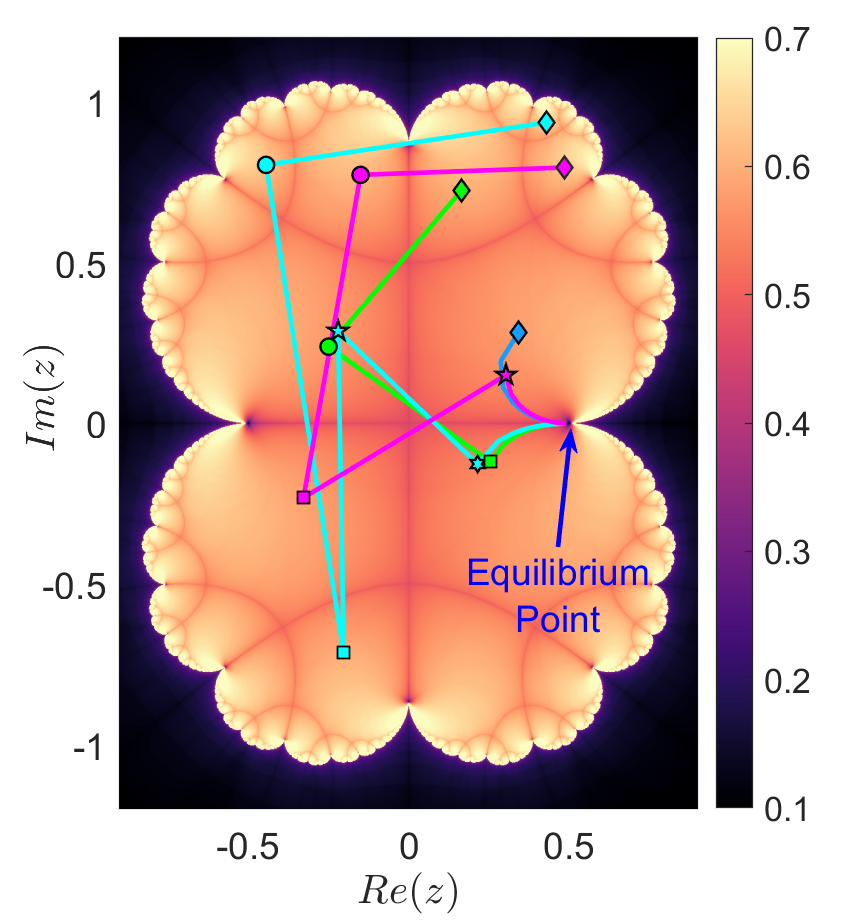}
		B)\includegraphics[scale=0.28]{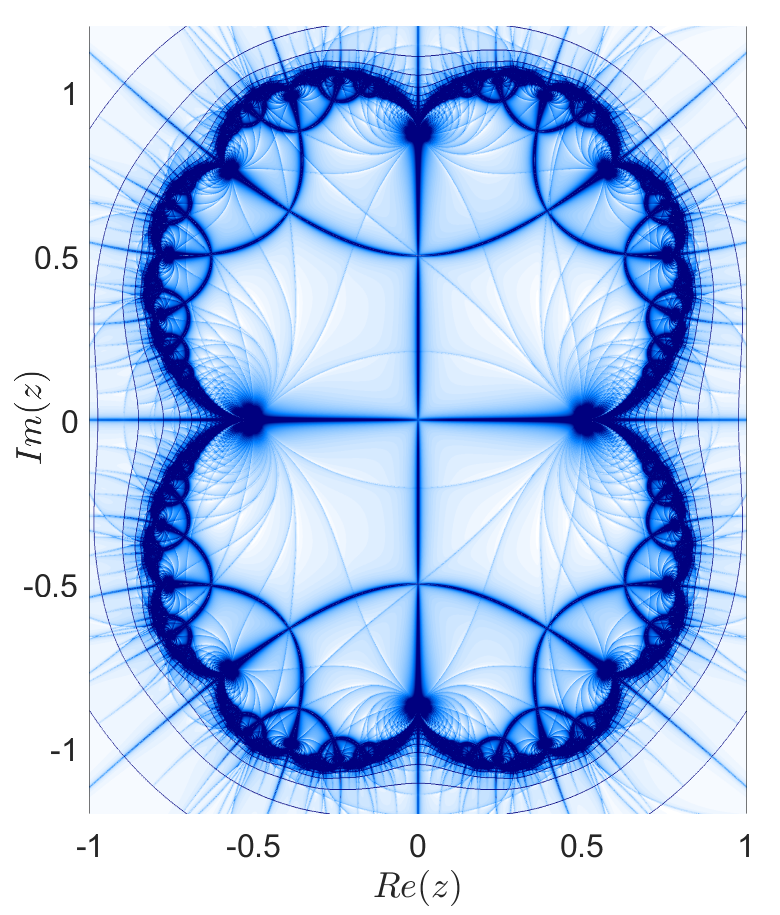}
	\end{center}
	\caption{Cauliflower Julia set ($c = 0.25$) revealed by DLDs using $p = 0.25$ and $N = 50$ iterations. A) Dynamics of several initial conditions, marked as diamonds. B) Tesselations, equipotential curves and external rays revealed by the ridges in the gradient of DLDs.}
	\label{LDs_cauli}
\end{figure}

\begin{figure}[!h]
	\begin{center}
		A)\includegraphics[scale=0.22]{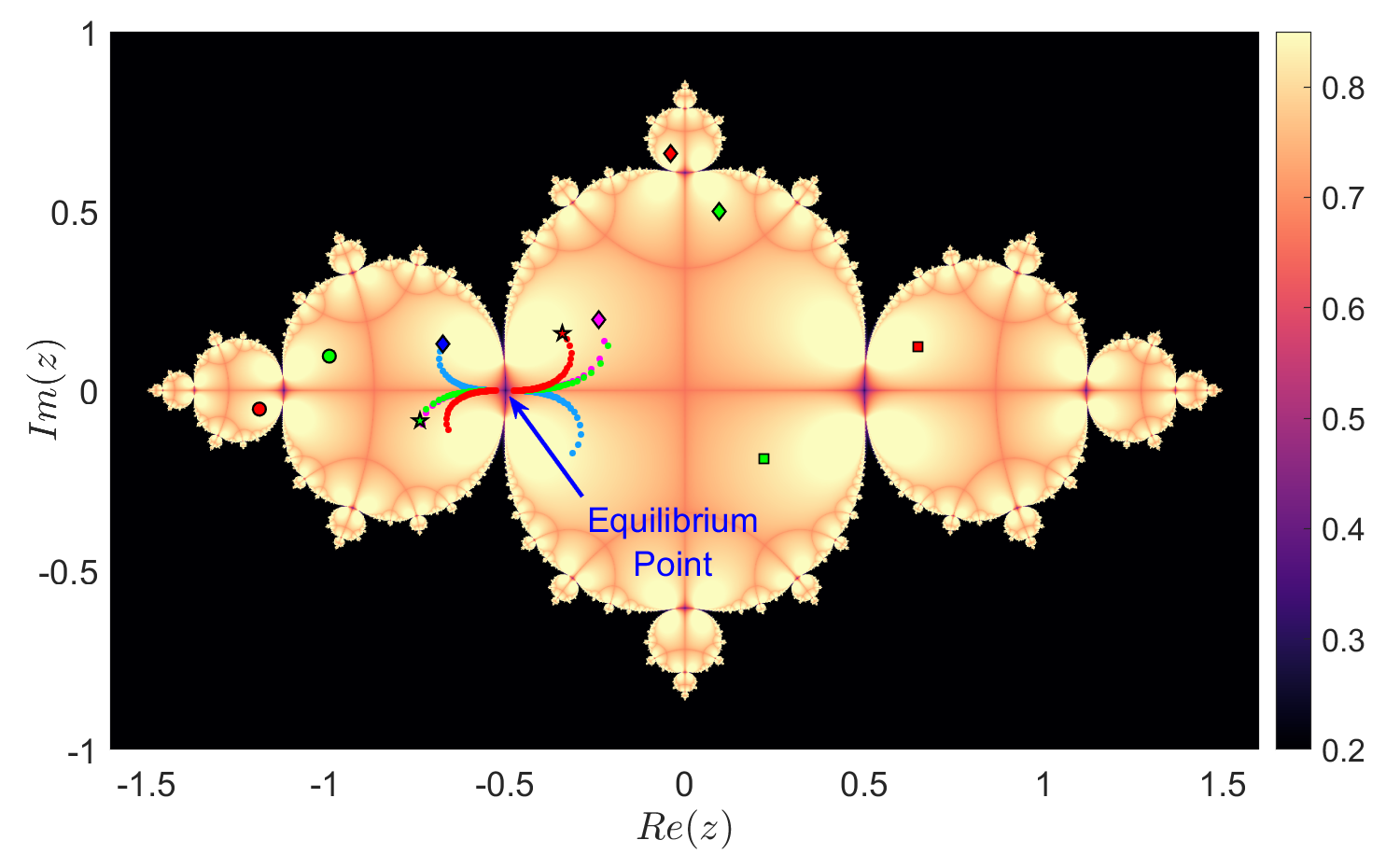}
		B)\includegraphics[scale=0.31]{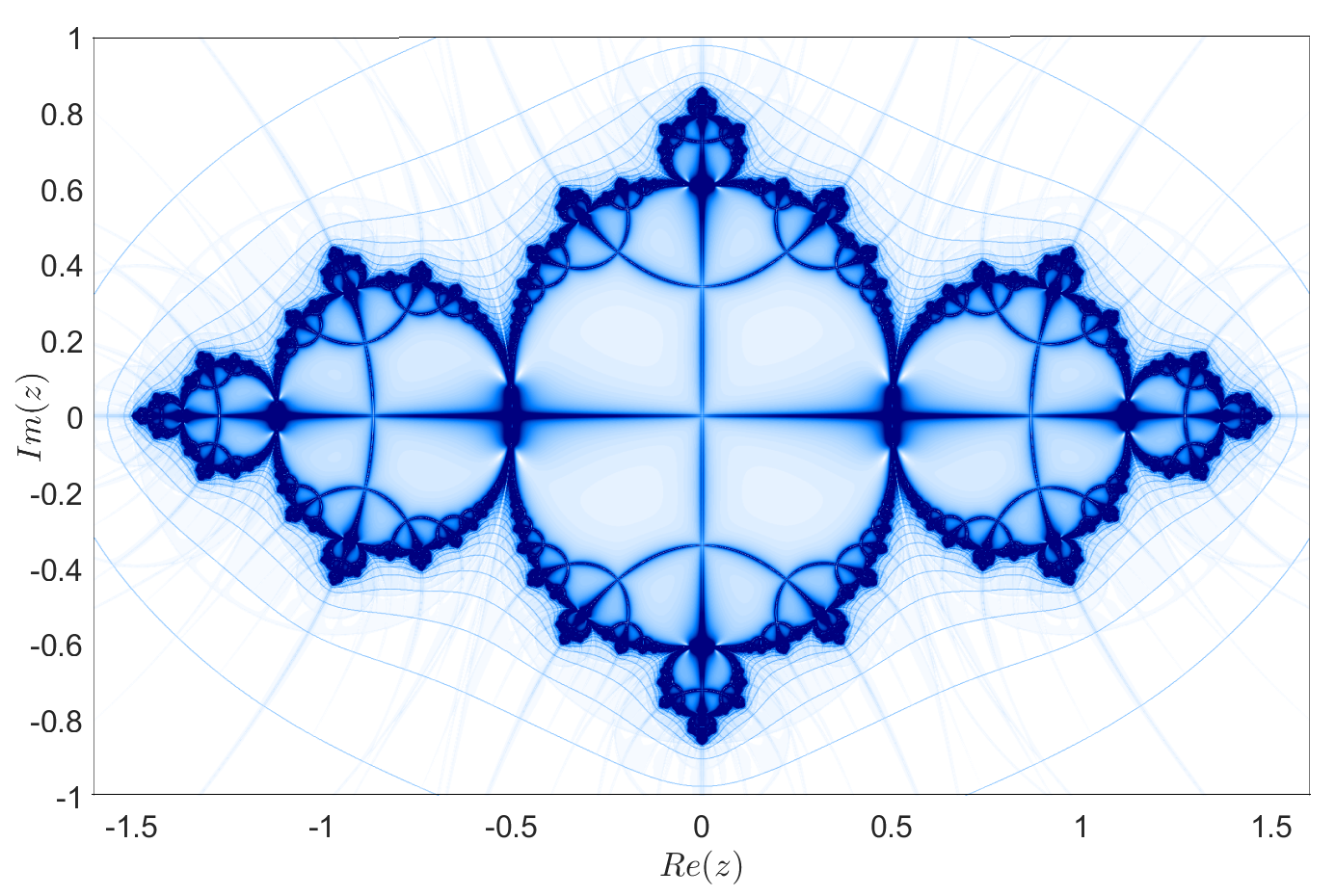}
	\end{center}
	\caption{Fat Basilica Julia set ($c = -0.75$) revealed by DLDs using $p = 0.25$ and $N = 100$ iterations. A) Dynamics of several initial conditions, marked as diamonds. B) Tesselations, equipotential curves and external rays revealed by the ridges in the gradient of DLDs.}
	\label{LDs_sanmarco}
\end{figure}

We turn our attention next to the analysis of the complex dynamics of parabolic Julia sets. Studies have revealed, see e.g. \cite{braverman} that standard programs used for the
computation of Julia sets work well when the underlying dynamical behavior in the system is hyperbolic, but experience high computational costs when dealing with parabolic dynamics due to an exponential-type slowdown in the neighborhood of parabilic fixed points. Furthermore, since some of the geometrical structures that characterize the internal dynamics of the set are very thin, they are very difficult to detect with precision. We demonstrate here that the method of DLDs described in this paper successfully circumvents these issues, and provides excellent results that make it a suitable tool to analyze parabolic-type dynamics efficiently. To illustrate this capabilites, we will take a look in particular to the Cauliflower ($c = 0.25$) and the fat Basilica ($c = -0.75$ Julia sets. We begin with the parameter value $c = 0.25$, and compute DLDs with an exponent $p = 1/4$ and $N = 50$ iterations. We depict in Fig. \ref{LDs_cauli} A) the scalar field obtained, and in panel B) the values where the gradient of DLDs takes very large values. Notice that besides the boundary of the fractal Julia set, we also successfully recover equipotential curves and external rays. Remarkably, and without further input, the method unveils the intricate tesselation pattern and parabolic sepals in the interior of the Cauliflower Julia set \cite{kawahira2003,kawahira2009a,kawahira2009b}. This geometrical template induces a checkerboard-type dynamics that we test by selecting initial conditions in different tiles of the underlying chkerboard. We mark initial conditions, $N = 0$, with diamonds and the successive iterations $N = 1,2,3$ with a circle, a square and a star respectively. In order to improve visualization and guide the eye, we have joined with straight lines the different iterations of a given orbit just to guide the eye, despite these lines do not have any mathematical meaning for the dynamics. note that the iterations jump from one tile to another with the same color, see the schemes included in \cite{kawahira2009a,kawahira2009b} until they reach the tile that is adjacent to the stable fixed point of the system, located at $z = 1/2$, and finally the orbit converges to the fixed point. A summary of this dynamical behavior is shown in Fig. \ref{LDs_cauli}. To provide further evidence, we also validate the numerical results that DLDs are giving by applying the same analysis discussed above to the case $c = -0.75$. A summary of the results obtained is dispayed in Fig. \ref{LDs_sanmarco}.

\begin{figure}[!h]
	\begin{center}
		A)\includegraphics[scale=0.215]{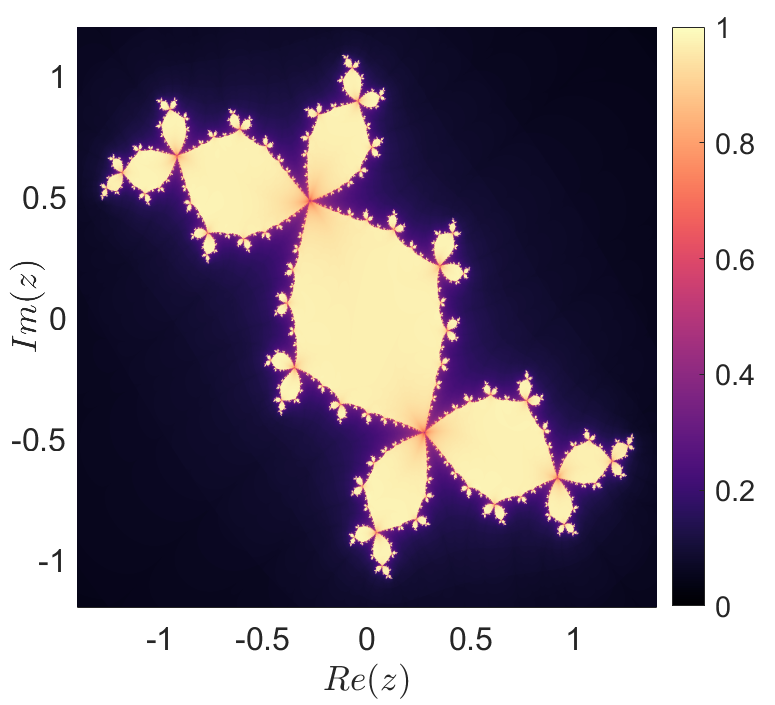}
		B)\includegraphics[scale=0.32]{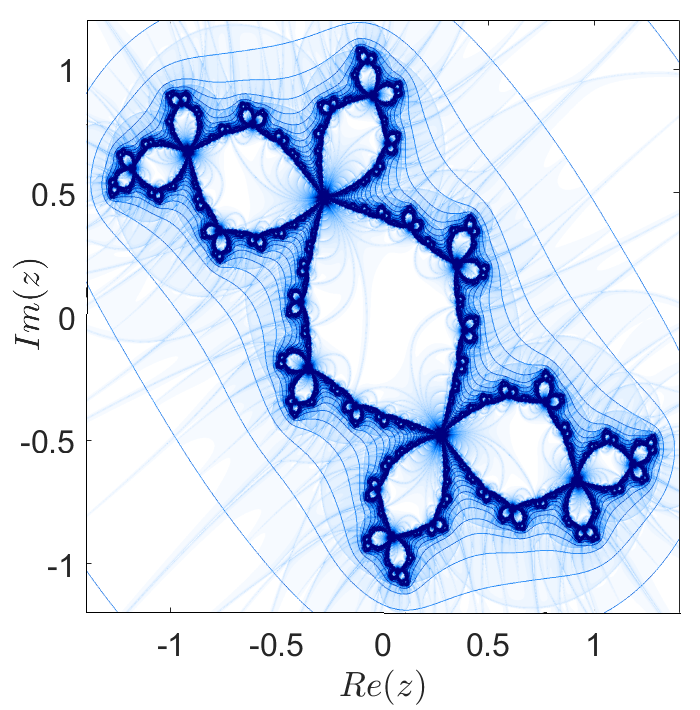}
		C)\includegraphics[scale=0.2]{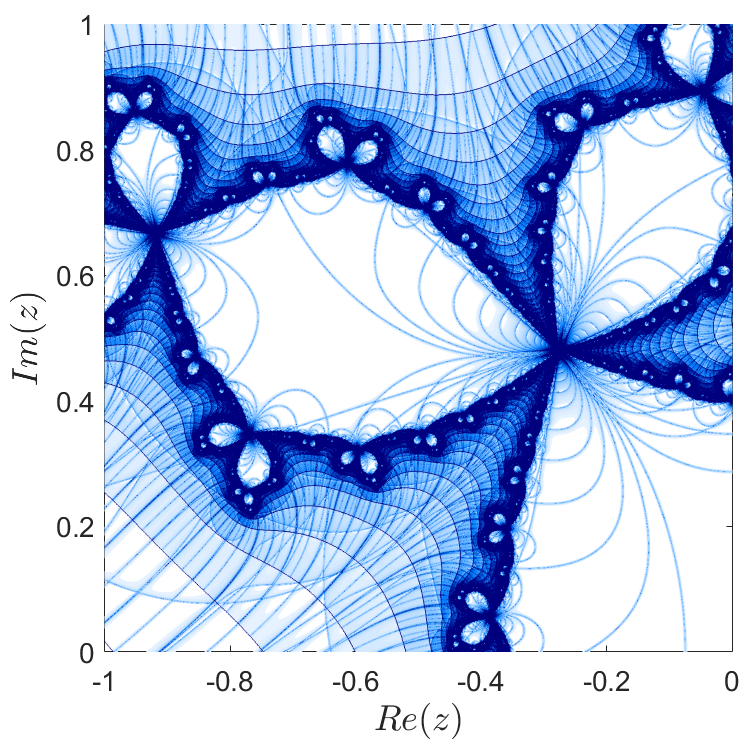}
	\end{center}
	\caption{Douady rabbit obtained by DLDs using $p = 0.25$ and $N = 50$ iterations for the quadratic mapping with $c = -0.123 + 0.745 \, i$. A) DLD scalar field; B) Dyanmical structures extracted from the gradient of DLDs, C) Zoom showing the intricate geometry of the internal laminations obtained in the Fatou components of the filled Julia set.}
	\label{LDs_douady}
\end{figure}

With the goal of highlighting the potential of the DLD scalar diagnostic to reveal Lyubich-Minsky laminations \cite{lyubich1997,kawahira2009b}, we apply it to the quadratic map with the well-known value $c = -0.123 + 0.745 \, i$ that generates the Douady rabbit fractal. In Fig. \ref{LDs_douady} we display the results obtained using DLDs with $p = 0.25$ and $N = 50$ iterations. Besides detecting, as in the previous cases considered, the external reys of the Julia set and the equipotentials, we can clearly see that the method also recovers the intricate geometry of the internal laminations in the Fatou components of the filled Julia set, see Fig. \ref{LDs_douady} C).

\begin{figure}[!h]
	\begin{center}
	    A)\includegraphics[scale=0.35]{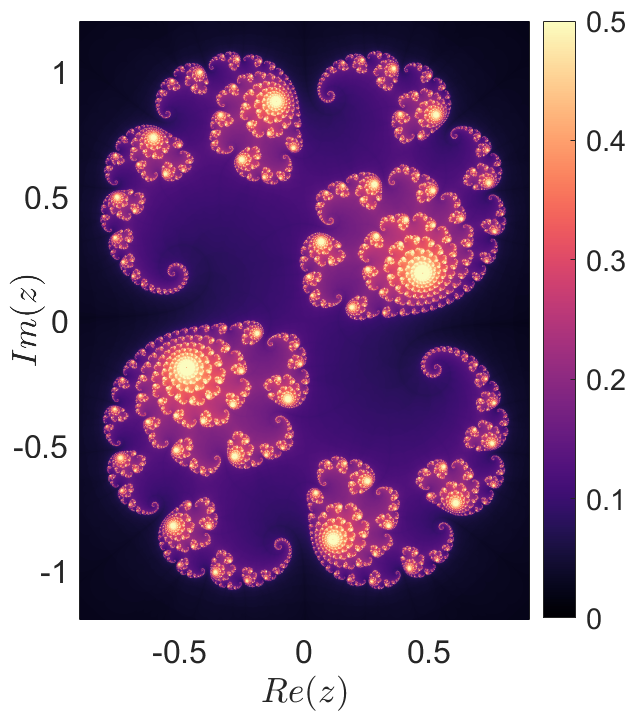}
		B)\includegraphics[scale=0.5]{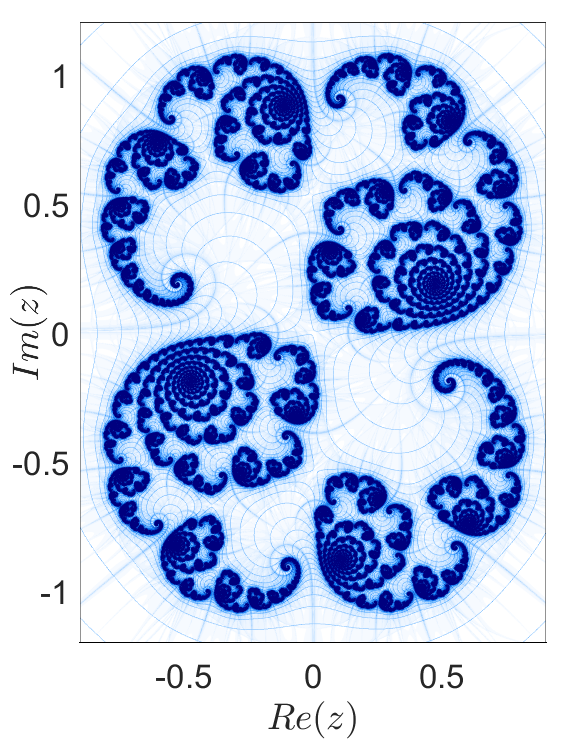} \\
		C)\includegraphics[scale=0.275]{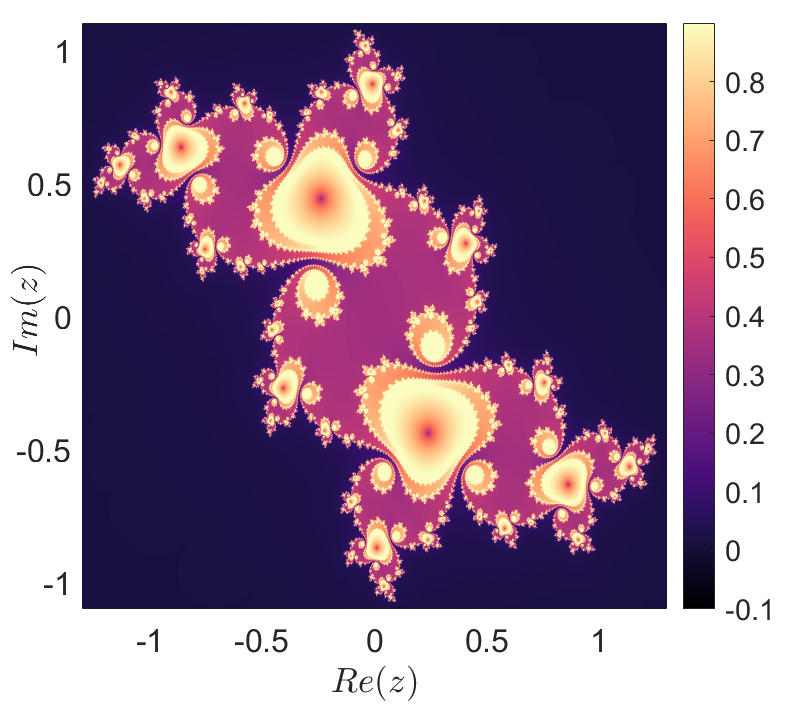}
		D)\includegraphics[scale=0.26]{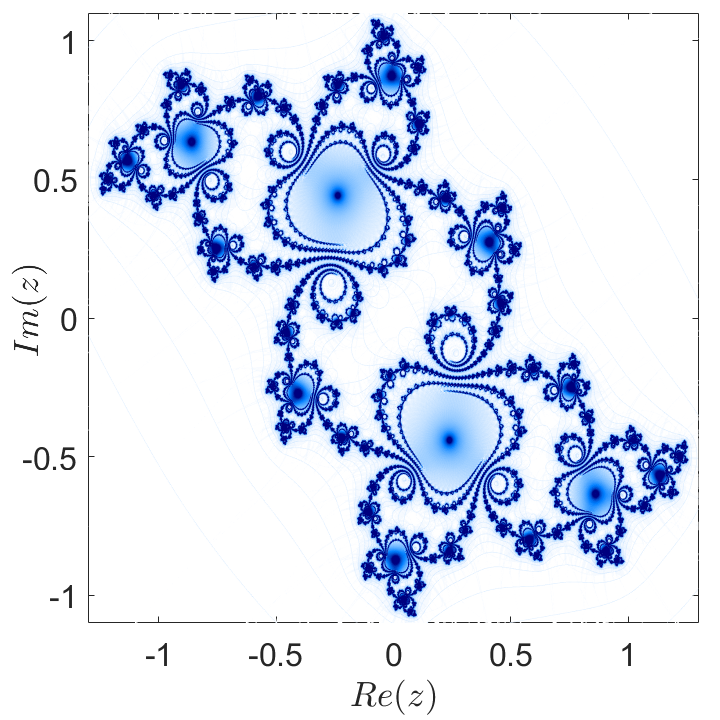} \\
		E)\includegraphics[scale=0.28]{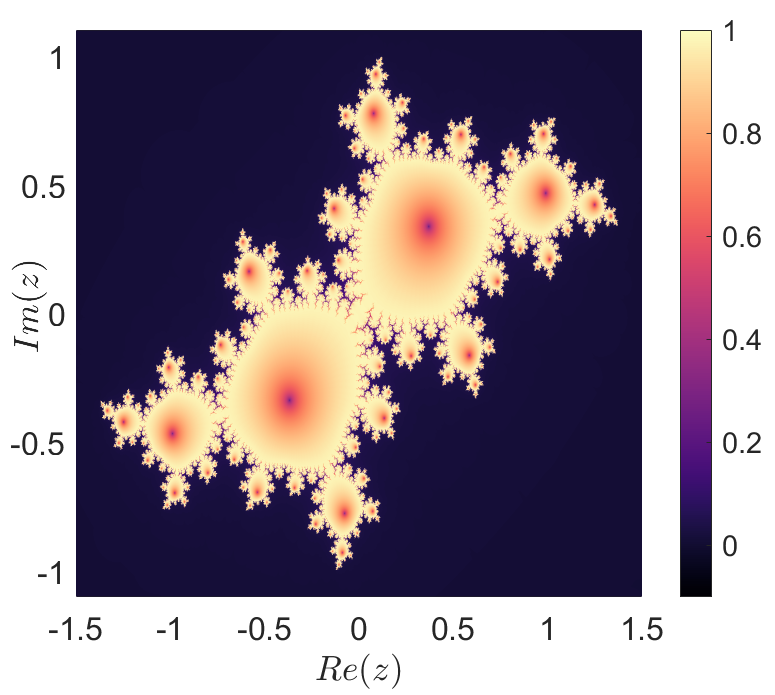}
		F)\includegraphics[scale=0.41]{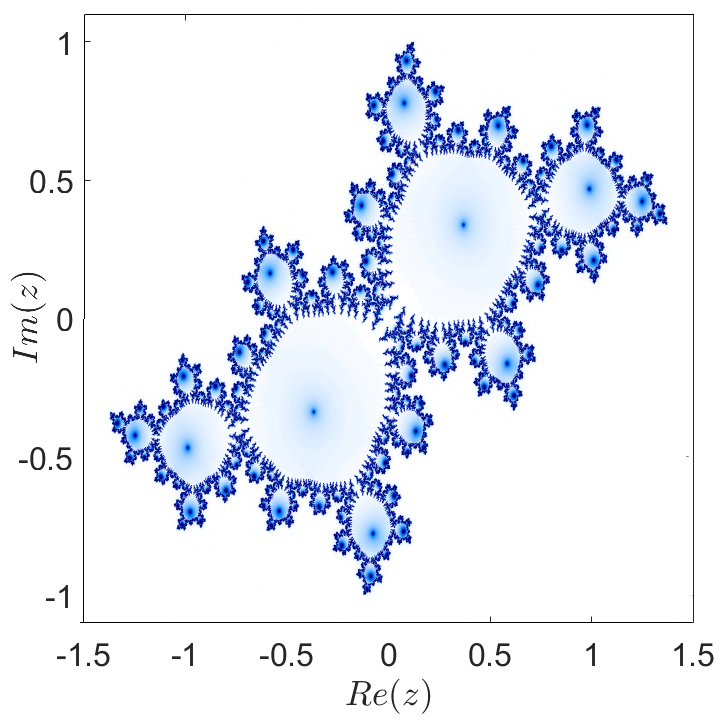} \\
	\end{center}
	\caption{In the left column we present DLDs using $p = 1/4$ and $N = 200$ iterations for the parameter values: A) $c = 0.285 + 0.01 \, i$ (open Cauliflower fractal); C) $c = -0.1 + 0.651 \, i$ (Flower-type fractal); E) $c = -0.391 - 0.587 \, i $ (Siegel-disk fractal). The right column depicts the corresponding fractal Julia sets extracted from the gradient of Lagrangian descriptors.}
	\label{LDs_fractals}
\end{figure}

We conclude the analysis of the dynamics induced by the qudratic map in Eq. \eqref{ds_qwad} by illustrating once more the capability of the diagnostic to unveil the feactal boundary of Julia sets at points where the scalar field it procduces is non-differentiable. We do so by reproducing the beautiful and intricate Julia sets given by the parameter values $c = 0.285 + 0.01 \, i$ (open Cauliflower fractal), $c = -0.1 + 0.651 \, i$ (Flower-type fractal) and $c = -0.391 - 0.587 \, i $ (Siegel-disk fractal). The results provided for this cases by applying DLDs using $p = 1/4$ and $N = 200$ iterations are displayed in Fig. \ref{LDs_fractals}. Observe the high-resolution that can be obtained for the intricate fractal boundary of the Julia sets with only a small number of iterations. It is important to remark that, as the number of iterations $N$ used to compute DLDs is increased, we are incorporating more information about the history of the orbits in the system, and therefore a richer geometrical description is achieved.

To finish this paper, we turn our attention now to a classical example where the complex dynamical system is defined by applying Newton's method to find the roots of a complex polynomial. In particular, we will show that LDs successfully reveal the classical Newton's fractal obtained when calculating the cubic roots of unity. Consider the polynomial $p(z) = z^3 - 1$ and write the discrete-time mapping:
\begin{equation}
z_{n+1} = f(z_n) = z_n - \dfrac{p(z_n)}{p^{\prime}(z_n)} = z_n - \dfrac{z_n^3 - 1}{3z_n^2} \quad,\quad n \in \mathbb{N} \cup \lbrace 0 \rbrace \;.
\label{newton_ds}
\end{equation}
This map is non-invertible and its dynamics nicely captures the highly-sensitive nature of the convergence of Newton's method to the different roots of the polynomial, depending on the location in the complex plane of the initial condition selected. This sensitive behavior of a dynamical system to initial conditions, which is one of the hallmarks of chaos, gives rise in the context of complex dynamics to fractal structures in the phase space of the system. When we apply DLDs to the complex mapping described in Eq. \eqref{newton_ds} using $p = 0.25$ and $N = 100$ iterations, we recover successfully the fractal set that governs the convergence of initial conditions to the different roots of unity. This interlaced three-fold and self-similar geometrical skeleton, which determines the boundaries of phase space regions that characterize the distinct dynamical fate of orbits,  offers an iconic testimony of the mathematical beauty of the subject of complex dynamics.

\begin{figure}[!h]
	\begin{center}
		A)\includegraphics[scale=0.28]{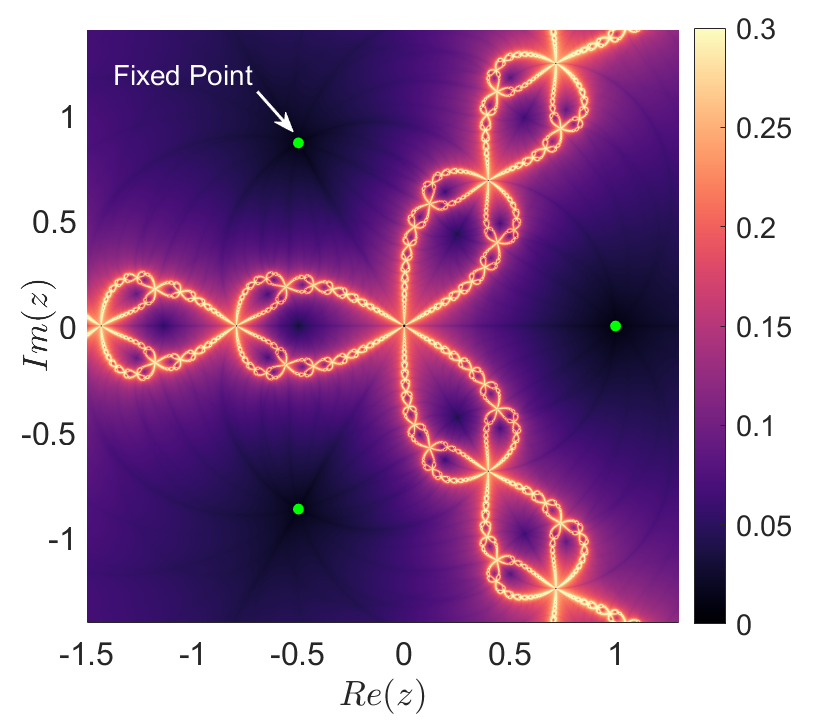}
		B)\includegraphics[scale=0.38]{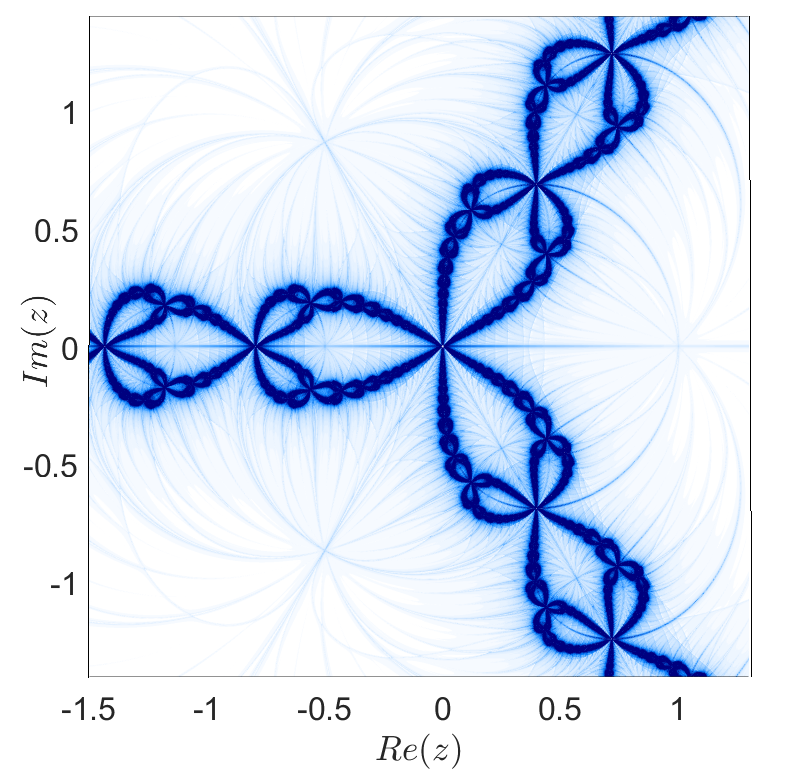}
	\end{center}
	\caption{A) Lagrangian descriptors calculated for the complex mapping defined in Eq. \eqref{newton_ds} using $p = 0.25$ and $N = 100$ iterations. We have marked the roots of unity, which are the fixed points of Newton's method, with green dots. B) Fractal boundary extracted from the gradient of Lagrangian descriptors.}
	\label{LDs_newton}
\end{figure}

\section{Conclusions}
\label{sec:conc}

In this work we have adapted the method of Discrete Lagrangian Descriptors so that it can become a useful scalar diagnostic to analyze the intricate phase space structures of complex dynamical systems. In order to demonstrate its capabilities, we have applied it to recover the fractal boundaries of Julia sets that arise from the classical quadratic mapping $z_{n+1} = z_n^2 + c$. Our analysis clearly shows that this tool proves to successfully unveil the boundaries between regions that display qualitatively distinct dynamical behavior. Moreover, we have also carried out further validations by using this methodology to study the dynamics generated by Newton's iterative method to calculate the roots of complex polynomials. All the results that we have discussed here provide us with good indicators that this methodology could become an extremely effective tool of the trade for the nonlinear dynamics community in their goal of exploring the relevant dynamical characteristics of complex mappings.

\section*{Acknowledgments} 

The author would like to acknowledge the financial support received from the EPSRC Grant No. EP/P021123/1 and the  Office of Naval Research Grant No. N00014-01-1-0769 for his research visits over the past two years to the School of Mathematics, University of Bristol. The inspiration for this work is the result of many fruitful discussions between the author and his collaborators Carlos Lopesino and Prof. Stephen Wiggins.

\bibliographystyle{natbib}
\bibliography{SNreac}

\end{document}